\documentclass[12pt]{article}
\usepackage{amsfonts}
\usepackage{latexsym}
\usepackage{color, amsmath,amssymb, amsfonts, amstext,amsthm, latexsym, dsfont}
\usepackage{epic,eepic}
\usepackage{graphicx}
\usepackage{epstopdf}
\usepackage{longtable}
\usepackage[title]{appendix}
\usepackage{indentfirst}
\usepackage{graphicx}
\usepackage{float}
\usepackage{graphicx}
\usepackage[utf8]{inputenc}
\usepackage{subfigure}
\usepackage{url}
\usepackage[colorlinks=false, linktocpage=true]{hyperref}
\usepackage{caption}
\allowdisplaybreaks[4]

\usepackage{tikz}
\usepackage{etoolbox}
\usepackage{enumitem}
\newcommand*{\circled}[1]{\lower.7ex\hbox{\tikz\draw (0pt, 0pt)%
    circle (.5em) node {\makebox[1em][c]{\small #1}};}}
\robustify{\circled}

\numberwithin{equation}{section}
 \numberwithin{Lem}{section}
 \numberwithin{Defi}{section}
 \numberwithin{Theo}{section}
 \numberwithin{Rem}{section}
  \numberwithin{Coro}{section}
  \numberwithin{Fig}{section}

\newcommand{\partialx}[1]{\frac{\partial}{\partial #1}}
\newcommand{\partialxy}[2]{\frac{\partial #1}{\partial #2}}

\begin{document}

\title{  A Data-Driven Approach for Discovering  the Most\\ Probable Transition Pathway for a\\ Stochastic Carbon Cycle System}

\author{\bf\normalsize{
Jianyu Chen$^{1,}$\footnotemark[1],
Jianyu Hu$^{1,}$\footnotemark[2],
Wei Wei$^{1,}$\footnotemark[3]
and Jinqiao Duan$^{2,}$\footnotemark[4]
}\\[10pt]
\footnotesize{$^1$Center for Mathematical Sciences, Huazhong University of Science and Technology,} \\
\footnotesize{Wuhan, Hubei 430074, China.} \\[5pt]
\footnotesize{$^2$Departments of Applied Mathematics \& Physics, Illinois Institute of Technology,}\\
\footnotesize{Chicago, IL 60616, USA.}
}

\footnotetext[1]{Email: \texttt{jianyuchen@hust.edu.cn}}
\footnotetext[2]{Email: \texttt{jianyuhu@hust.edu.cn}}
\footnotetext[4]{Email: \texttt{duan@iit.edu}}
\footnotetext[3]{Email: \texttt{weiw16@hust.edu.cn}}
\footnotetext[1]{is the corresponding author}

\date{}
\maketitle
\vspace{-0.3in}

\maketitle
\begin{abstract}
 Many natural systems exhibit tipping points where changing environmental conditions spark a sudden shift to a new and sometimes quite different state. Global climate change is often associated with the stability of marine carbon stocks. 
 We consider a stochastic carbonate system of the upper ocean to capture such transition phenomena. Based on the Onsager-Machlup action functional theory, we calculate the most probable transition pathway between the metastable and oscillatory states via a  neural shooting method, and further explore the effects of external random carbon input rates on the most probable transition pathway, which provides a basis to recognize naturally occurring tipping points. Particularly, we investigate the effect of the transition time on the transition pathway and  further compute the optimal transition time using physics informed neural network, towards the maximum carbonate concentration state in the oscillatory regimes. This work offers some insights on the effects of random carbon input on climate transition in a simple model.


\end{abstract}

{\bf Lead paragraph:}

\textbf{Climate nowadays has been a major concern because of the significant impact of carbon dioxide from human activities. The oceanic carbonate system plays an important role in the global carbon cycle. In this study, we consider a stochastic oceanic carbonate system, which has a metastable state and oscillation regimes. Due to the random fluctuations in the external carbon input rate, the system undergoes a transition between the two regimes. We apply a neural shooting method to compute the most probable transition pathway from the metastable state to the oscillation regimes, based on the Onsager-Machlup action functional theory. Moreover, we compare the effects of different carbon input rates on the transition mechanism. The results show that not only dissolved inorganic carbon (DIC) $w$, but also $[CO^{-2}_{3}]$ concentration $c$ is fluctuated with external carbon dioxide input rate $\nu$, and there has been a huge oscillation at $\nu$=0.1, which can be considered as a critical transition. When $\nu$ is growing from 0.2, the value of $c$ fluctuates around $58$  $\mu mol\cdot kg^{-1}$ without a significant abrupt change. Moreover, we detect the impacts of transition time $T$ on transition pathway when $\nu$=0, which indicates that as time increases, the endpoint of the most probable transition pathway moves toward the right half of the limit cycle. That is, for the stochastic carbon cycle system, with the increase of time, the system is more and more inclined to transition from the metastable state to an oscillatory state with increasing $[CO^{-2}_{3}]$ concentration $c$. Furthermore, we explore the optimal transition time for a typical transition phenomenon, and it plays an important role on the early warning signal forecast for disruptions. This is achieved by a physics informed neural network.}
\medskip

\section{Introduction}
The Earth's carbon cycle system is a link between photosynthesis and respiration. But the occurrence of emergencies breaks this calm and stable cycle, which is shown in geological records as relatively sudden and large changes in the composition of sedimentary carbon isotopes \cite{1999Interpreting}. Such events are usually attributable to changes in carbon fluxes and concentration \cite{1999Interpreting,Rothman2003Dynamics}.

The extent to which the carbon cycle system is destroyed is determined not by the intensity of the cycle disturbance, but by the internal dynamics of the system itself. Once the carbon dioxide added to the oceans exceeds the threshold, the environment will change dramatically, and this catastrophic change will bring about sudden changes in the human living environment.

But the cause of the destruction of the Earth's environment remains a mystery. Each change will lead to a temporary change in the ocean's carbon reserves. Although the causes are controversial, these changes are often interpreted as a proportional response to external carbon dioxide inputs. Rothman \cite{Rothman14813} argues that the magnitude of many of the damage is determined not by the strength of external disruption, but by the internal dynamics of the carbon cycle. Therefore, it is of great significance to study the effect of the change of carbon flux entering the ocean from outside on the system.

Rothman \cite{Rothman14813} formulated a deterministic mathematical model to investigate the dynamical behaviors of the carbon cycle system. 
The disruptions were suggested to be initiated by perturbation of a permanently stable steady state beyond a threshold. However, he did not take random fluctuations in the strength of an external source of $CO_2$, as in volcanic emissions and human activities, into account. While in \cite{wei2021most}, a quite small random perturbation was considered to describe the noisy fluctuations, which is the large deviation limit regime, where infinitesimal perturbations are considered. More often, we must take into account the general random fluctuations rather than infinitesimal random perturbations. In this sense, we consider a stochastic carbon cycle model, which includes inevitable general random fluctuations from the external carbon input.

We look into the bistable case of the stochastic carbon cycle system. When the bifurcation parameter $c_x$ is in the interval $ (56, 62.61)$ in Rothman's case, it has one stable steady state and one stable limit cycle, with a sandwiched unstable limit cycle.
Random fluctuations result in transition phenomena of this type of system, which are impossible in the deterministic case. We are interested in transition behaviors from the metastable state to the oscillatory regimes for the stochastic carbon cycle system. 
The Onsager-Machlup theory is an effective tool to capture transition behaviors of stochastic dynamical systems, which has been widely investigated \cite{capitaine1995onsager,fujita1982onsager,ikeda2014stochastic,shepp1992note,takahashi1981probability} . 
More precisely, the Onsager-Machlup action functional gives an estimation of the probability of the solution paths on a small tube. The minimum value of the action functional means the largest probability of the solution paths on the tube. By minimizing the action functional, we can obtain the most probable transition pathway connecting two states, which satisfies the Euler-Lagrange equation.

We propose the neural shooting method \cite{Hu2021TransitionPF,ibraheem2011shooting} to numerically calculate the Euler-Lagrange equation and obtain the most probable transition pathway. Moreover, we investigate the effect of time $T$ on the transition pathway and compute the optimal transition time using physics informed neural network (PINNs) \cite{Hu2021DatadrivenMT,raissi2019physics} for a fixed state on the oscillatory regimes in the stochastic carbon cycle system, which is of great significance for practical earning warning forecasting. The shooting method is an effective method to solve the transition between two states in fixed time. Because it numerically solves the difficulty of the two-point boundary value problem.

This paper is organized as follows. In section 2, we consider a stochastic model to investigate the dynamical behaviors of the carbon cycle system. An Onsager-Machlup approach is presented to study transition from a metastable state to an oscillatory regime. In section 3, we propose both the neural shooting method and PINNs method to calculate the most probable transition pathway. Particular attention is given to the study of the effect of time on the most probable transition pathway. We then discuss the optimal transition time for a fixed state of the stochastic carbon cycle system. The conclusion highlights strengths and weakness of our research in section 4.

\section{Stochastic carbon cycle system}
 

\subsection{A marine carbonate system }
In Rothman \cite{Rothman14813}, the dissolved inorganic carbon (DIC) $w$ and carbonate ions $CO_3^{2-}$ are concentrations. He introduced a deterministic model to describe dynamical behaviors of the carbonate system as follows.
\begin{equation}
\begin{aligned}
dc&=\left[\mu\left[1-b s\left(c, c_{p}\right)-\theta \bar{s}\left(c, c_{x}\right)-\nu \right]+w-w_{0}\right]f(c)dt,\\
dw&=\left[\mu\left[1-b s\left(c, c_{p}\right)+\theta \bar{s}\left(c, c_{x}\right)+\nu\right]-w+w_{0}\right]dt,
\end{aligned}
\end{equation}
where $\mu$ is a characteristic concentration, $b$ is the maximum $CaCO_3$ burial rate, $\theta$ is the maximum respiration feedback rate, $\nu$ is the injection rate of $CO_2$, $w$ is the reference DIC concentration. Here, $\bar{s}=1-s$, $s$ is the sigmoid function, and $f$ is the “buffer function”
\begin{equation}
s\left(c, c_{p}\right)=\frac{c^{\gamma}}{c^{\gamma}+c_{p}^{\gamma}},\ \  f\left(c\right)=f_0\frac{c^{\beta}}{c^{\beta}+c_{f}^{\beta}},
\end{equation}
where $c_x$ is the crossover $[CO^{-2}_{3}]$ (respiration), $c_p$ is the crossover $[CO^{-2}_{3}]$ (burial), $c_f$ is the crossover $[CO^{-2}_{3}]$ (buffering), and $f_0$ is the maximum buffer factor and $\gamma,\beta$ are the sigmoid sharpness indexes.

However, Rothman did not take random fluctuations in  the strength of an external source of $CO_2$, as in volcanic emissions and human activities into account. As in \cite{wei2021most}, a small random perturbation was considered to describe the noisy fluctuations. More often, we must take into account the general random fluctuations rather than infinitesimal random perturbations. To this end, we consider the following stochastic carbon cycle system model
\begin{equation}
\begin{aligned}
dc&=\left[\mu\left[1-b s\left(c, c_{p}\right)-\theta \bar{s}\left(c, c_{x}\right)-\nu \right]+w-w_{0}\right]f(c)dt-\mu f(c)dB_t^1,\\
dw&=\left[\mu\left[1-b s\left(c, c_{p}\right)+\theta \bar{s}\left(c, c_{x}\right)+\nu\right]-w+w_{0}\right]dt+\mu dB_t^2.
\end{aligned}\label{carbon}
\end{equation}


In this stochatic carbon cycle system, time unit is $10^4$ years (so t=1, means the time is 10000 years), and $c, w$ have unit $\mu mol\cdot kg^{-1}$.

The parameters $\mu$, $b$, $\theta$, $c_p$, $c_x$, $c_f$, $w_0$, $\gamma$ and $\beta$ are constants and are set to fit in properties of the modern ocean.The value of these parameters are listed in \cite[SI Appendix, Table S1]{Rothman14813}.  We only change the value of $c_x$ and $\nu$ in this article.

\begin{figure}[htbp]
    \centering
    \includegraphics[scale=0.2]{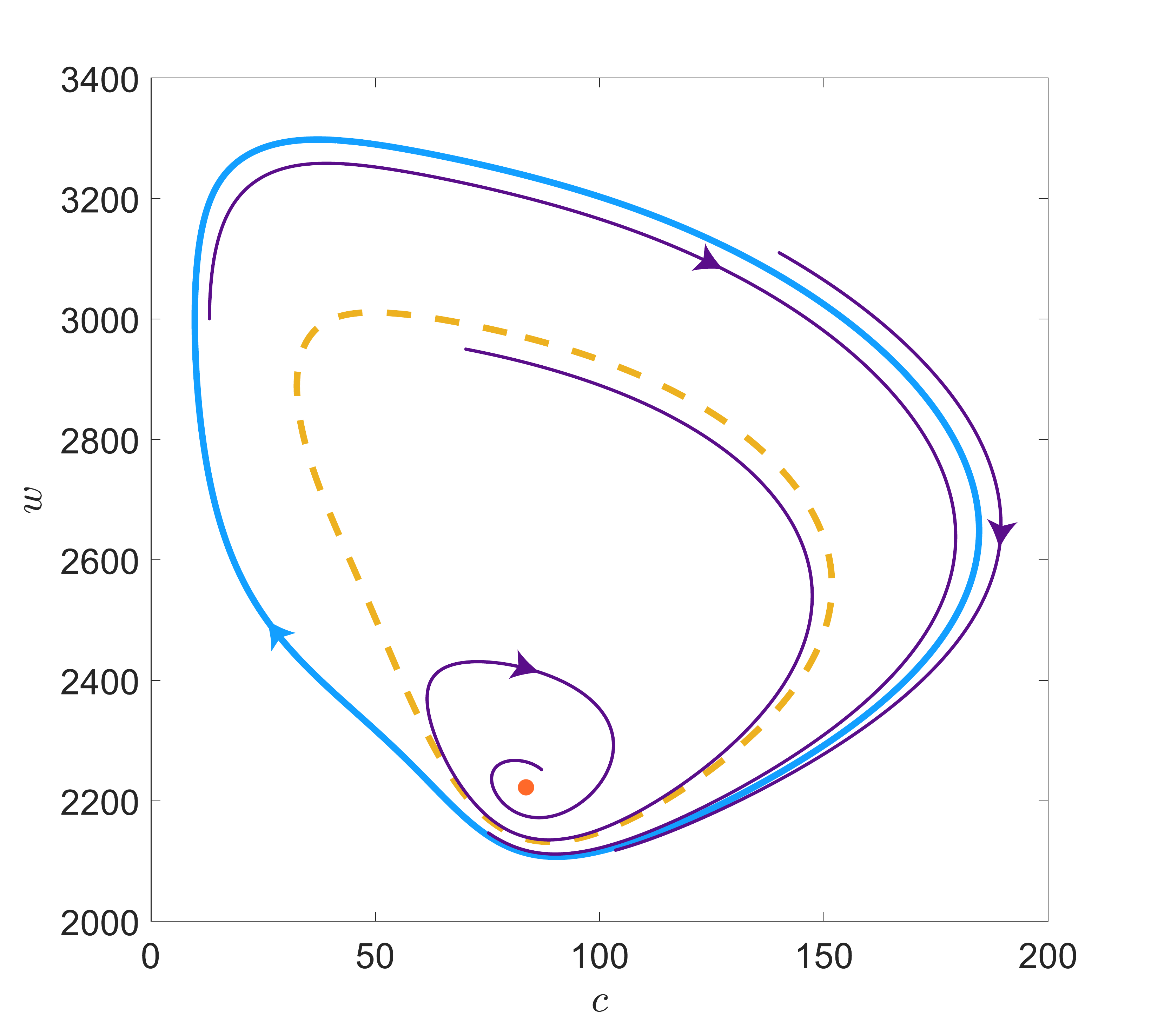}
    \caption{Phase-space trajectories in the bistable regime. Here $c_x=58 \mathrm{ \mu mol \cdot kg^{-1}}$. The yellow dashed limit cycle is unstable. Trajectories (purple arrow lines) initialized inside the unstable limit cycle return to the stable fixed point (the reddish orange point). Trajectories (purple arrow lines) initialized outside the unstable limit cycle evolve to the stable limit cycle (the blue arrow loop).}
    \label{doublecycle}
\end{figure}

In the bistable regime $c_x \in (56, 62.61)$, there are one stable state, one stable limit cycle and one unstable limit cycle. We take $c_x=58$. Phase-space trajectories in the bistable regime are shown in Fig. \ref{doublecycle} \cite{wei2021most}. The input rates $\nu$ of $CO_2$ are taken as $0, 0.1, ..., 0.9$.



 


\subsection{An Onsager-Machlup approach}
In this section, we briefly introduce the Onsager-Machlup theorem and the most probable transition pathway.

We consider a generalized stochastic differential equation in $\mathbb{R}^d$:
\begin{equation}
\begin{split}\label{SDE}
dX(t)=\tilde{f}(X(t))dt+\sigma(X(t))dB(t), t\in[0,1],
\end{split}
\end{equation}
with initial data $X(0)=x_0 \in \mathbb{R}^d$, where $\tilde{f}:\mathbb{R}^d\rightarrow\mathbb{R}^d$ is a regular function, $\sigma:\mathbb{R}^d\rightarrow\mathbb{R}^{d\times k}$ is a $d\times k$ matrix-valued function, and $B$ is a Brownian motion in $\mathbb{R}^k$. For the stochastic carbon cycle system \eqref{carbon}, $\sigma(x,y)=
\left(
\begin{matrix}
   -\mu f(x)  &  0\\
   0  & \mu
\end{matrix}
\right)
$, and $B=(B^1,B^2)$.

The well-posedness of stochastic differential equation \eqref{SDE} has been widely investigated \cite{karatzas2012brownian}. We are interested in the transition phenomena between two metastable states (taking as the stable equilibrium states of the corresponding deterministic system). 

To investigate transition phenomena, one should estimate the probability of the solution paths on a small tube. The Onsager-Machlup theory of stochastic dynamical system \eqref{SDE} gives an approximation of the probability
\begin{equation}
\begin{split}
\mathbb{P}(\{\|X-z\|_T \leqslant \delta\}) \propto C(\delta,T) \exp \left\{-S(z,\dot{z})\right\},
\end{split}
\end{equation}
where $\delta$ is positive and sufficiently small, $C(\delta,T)=\mathbb{P}(\{\|B\|_T \leqslant \delta\})$, the probability of Bownian motion in the $\delta$-tube, $\|\cdot\|_T$ is the uniform norm of the space of all continuous functions in the time interval $[0,T]$, $z$ is function in this space and the Onsager-Machlup action functional is
\begin{equation}
\begin{split}\label{action}
S(z,\dot{z})=\frac{1}{2}\int_0^T[(\dot{z}-b(z))V(z)(\dot{z}-b(z))+\operatorname{div}b(z)-\frac{1}{6}R(z)]dt,
\end{split}
\end{equation}
where $V(z)=(\sigma(z)\sigma^{*}(z))^{-1}$, $R(z)$ is the scalar curvature with respect to the Riemannian metric induced by $V(z)$, $b^i(z)=\tilde{b}^i(z)-\frac{1}{2}\sum\limits_{l,j}(V^{-1}(z))^{lj}\Gamma_{lj}^i$ is the $i$ component of $b$. Here, $\Gamma_{lj}^i$ is the Christoffel symbols associated with this Riemannian metric, which satisfies 
\begin{equation}
\begin{split}
\Gamma_{lj}^i=\frac{1}{2}\sum\limits_m g^{im}(\partialx{x^j}g_{lm}+\partialx{x^l}g_{jm}-\partialx{x^m}g_{lj}),\label{eq:2.7}
\end{split}
\end{equation}
where $(g^{ij})(z)$ is the inverse of the Riemannian metric $(g_{ij})(z)=V(z)$. The divergence $\operatorname{div}b(z)$ is defined as
\begin{equation}
\begin{split}
\operatorname{div}b(z)=\frac{1}{\sqrt{|V(z)|}}\sum\limits_i\partialx{z_i}\left(b^{i}(z)\sqrt{|V(z)|}\right),
\end{split}
\end{equation}
where $|V(z)|$ is the determinate of Riemannian metric. The Onsager-Machlup action functional can be considered as the integral of a Lagrangian 
\begin{equation}
\begin{split}
L(z,\dot{z})=(\dot{z}-b(z))V(z)(\dot{z}-b(z))+\operatorname{div}b(z)-\frac{1}{6}R(z).
\end{split}
\end{equation}
For the stochastic carbon cycle system \eqref{carbon}, the Onsager-Machlup action functional function

The minimum value of the action functional \eqref{action} means the largest probability of the solution paths on a small tube. 
We can investigate the most probable transition pathway by minimizing the the Onsager-Machlup action functional.
It captures sample paths of the largest probability around its neighbourhood.
In the carbon cycle system, we consider the transition from a metastable state to an oscillatory regime. The metastable state is the stable state $z^{*}$ and the oscillatory regime is a limit cycle $M$ in the corresponding deterministic differential equation. Thus , we need to find the minimizer of the following optimization problem
\begin{equation}
\begin{split}\label{mini}
\inf_{y \in M}\inf_{z \in \bar{C}_{z^{*}}^{y}(0,T)} S(z,\dot{z}),
\end{split}
\end{equation}
where $\bar{C}_{z^{*}}^{y}(0,T)$ is the space of all absolutely continuous functions that start at $z^{*}$ and end at $y$.

To get the minimizer of the optimization problem \eqref{mini}, we will take the infimum step by step. For each $y\in M$, by the variational principle, the most probable transition pathway from $z^{*}$ to $y$ satisfies the Euler-Lagrange equation
\begin{equation}
\begin{split}\label{EL_eqn}
\frac{d}{dt}\partialx{\dot{z}}L(z,\dot{z})=\partialx{z}L(z,\dot{z}),
\end{split}
\end{equation}
with initial state $z(0)=z^{*}$ and final state $z(T)=y$. The detailed derivation of the Euler-Lagrangian equation is in Appendix \ref{AppA} for a more close consultation. Comparing the action functional values among all $y\in M$, we obtain the most probable transition pathway from a metastable state to an oscillatory regime. In next section, we will discretize the limit cycle and propose both the neural shooting method and PINNs method to numerically calculate the Euler-Lagrange equation \eqref{EL_eqn} to explore the stochastic carbon cycle system.


\section{ A data-driven approach for simulating and computing}
In this section, we present both neural network shooting method and PINNs method to explore the stochastic carbon cycle system. The shooting method is used for computing the most probable transition pathway and reflecting a phenomenon about the effect of time on the transition. The PINNs method do forecast the optimal transition time.

\subsection{A neural shooting method}
The Euler-Lagrange equation \eqref{EL_eqn} for our stochastic  carbon cycle system is actually a second order ordinary differential equation of two dimensions.
\begin{equation}
\begin{split}\label{EL_eqn1}
\ddot{z}=g(z,\dot{z}),
\end{split}
\end{equation}
where $g$ is associated to the drift term $f$ and diffusion term $\sigma$ in stochastic differential equation \eqref{SDE}.

Combining the initial point and the finial point, it turns to a two-point boundary value problem, which equals to a Cauchy problem in four dimensions
\begin{equation}
\begin{split}\label{Cauchy}
\dot{x}&=y,\\
\dot{y}&=g(x,y),
\end{split}
\end{equation}

However, for this type of Cauchy problem, we do not know its initial velocity condition, but rather the end position condition.

Shooting method is an effective tool to solve the two-point boundary value problem. This paper intends to adopt a neural network approach to handle the tricky problem. 
The main idea is to construct a neural network and approximate the mapping from the final point to the initial velocity, such that the solution of the Cauchy problem \ref{Cauchy} with this initial velocity at time $T$ meets the target given final point.
Therefore, we need to generate a data set to supply the training basis for neural network.

\noindent \textbf{Part 1 - Data Generation}

\begin{figure}[htbp]
\begin{minipage}[]{0.5 \textwidth}
 \leftline{~~~~~~~\tiny\textbf{(\ref{fig:2}a)}}
\centerline{\includegraphics[width=6.5cm]{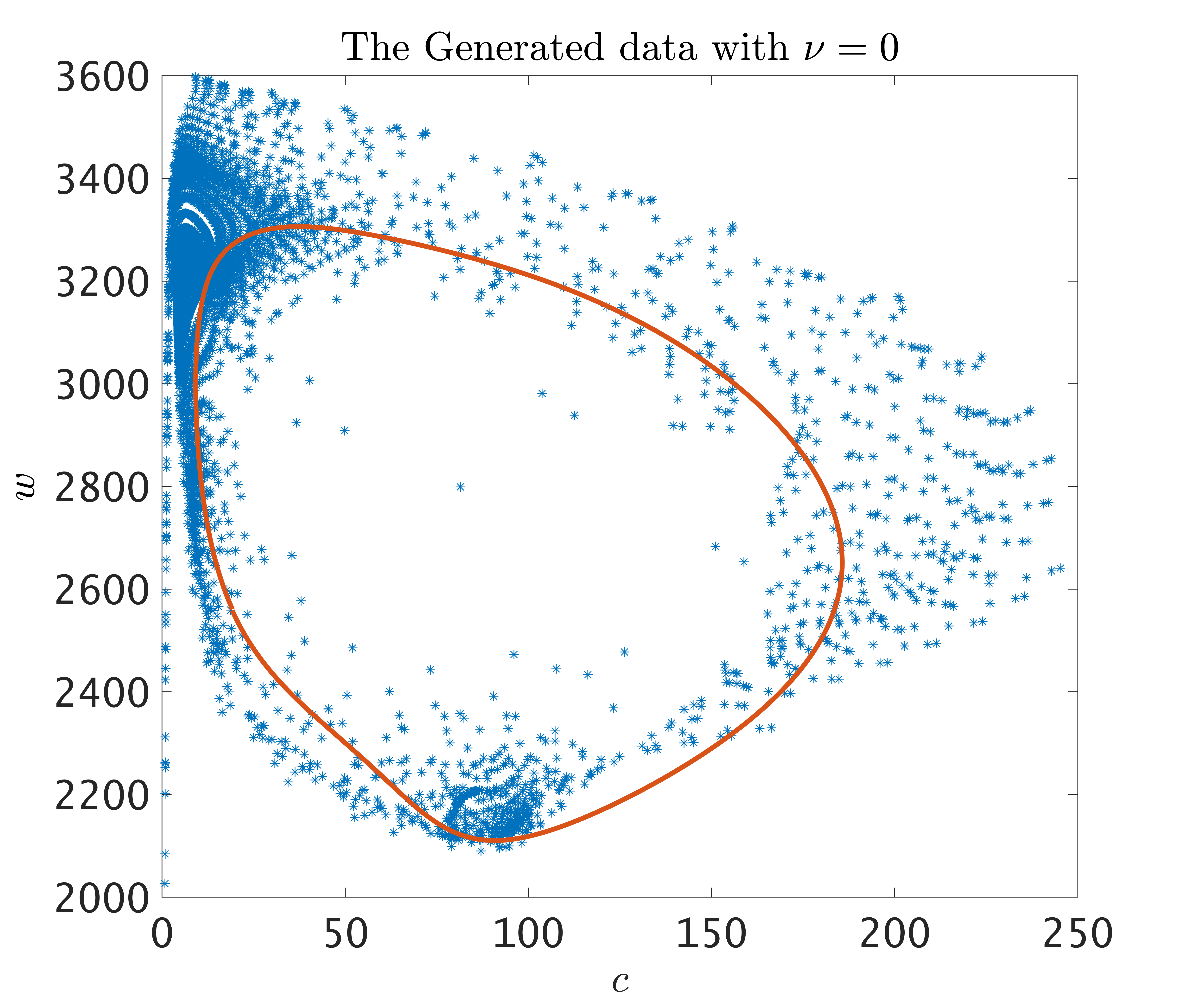}}
\end{minipage}
\hfill
\begin{minipage}[]{0.5 \textwidth}
 \leftline{~~~~~~~\tiny\textbf{(\ref{fig:2}b)}}
\centerline{\includegraphics[width=6.5cm]{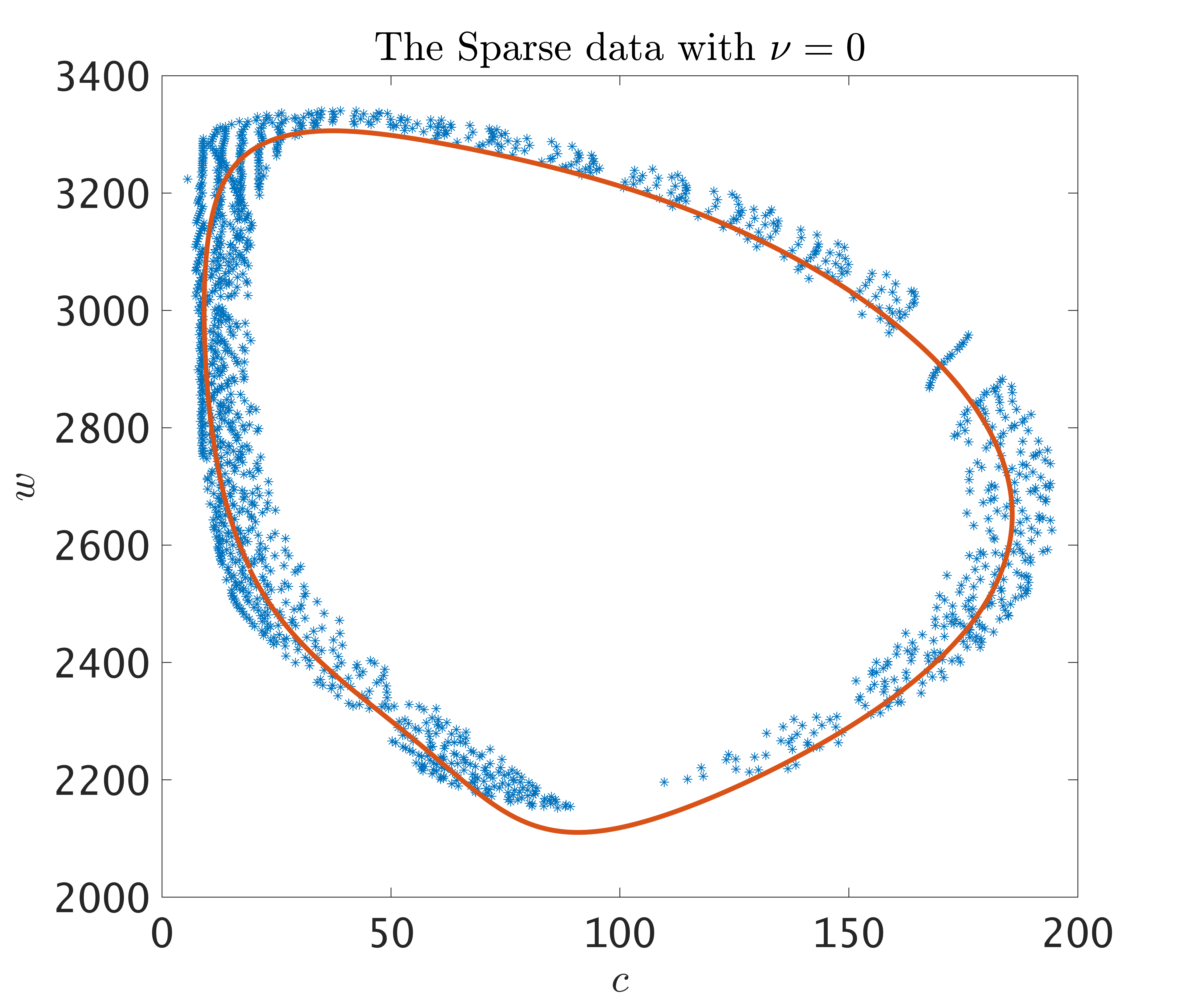}}
\end{minipage}
\caption{The data is generated by solving the Cauchy problem via discrete Euler method when the input rate $\nu$=0. (2a): the original data and the right-hand side is sparsely processed. The blue points stand for the final points. (2b): the limit cycle of the deterministic carbon cycle system.}\label{fig:2}

\end{figure}

We use the Euler method to solve the Cauchy problem with a given set of initial velocities. Thus, we obtain a set of final points at time $T$ which corresponds to a set of initial velocities.  This dataset consists of final point-velocity pairs. We try to find out such a mapping structure from the dataset through a neural network and then  predict the initial velocities corresponding to the limit cycle sets. Obviously, we prefer to choose datasets where the target points locate near the limit cycle, which makes the learning of the neural network more efficient.

Fix the parameter $\nu$ and generate the data by Euler method. For example, when $\nu = 0$, the generated data is shown as Fig. \ref{fig:2}. For $\nu$ of other values, the data's distributions are similar to the distribution shown as Fig. \ref{fig:2} and are not listed here. Readers who are interested in it can refer to the links in the data availability for more details. As can be seen from the distribution in Fig. \ref{fig:2}, we will use the blue scattered points distributed around the limit cycle as the input of the neural network and the corresponding initial velocities as the output to establish a mapping.

\noindent \textbf{Part 2 - Neural Shooting}

We adopt the fully connected neural network`` 2-8-16-8-2'' to approximate the mapping. We train three thousand times on each data set and finally use a model with a generalization error of $0.1$ to make the prediction and compute the most probable transition pathway. 

 For the reachability of the transition pathway, we give the following explanation. For a given point $(x,y)\in M$, we predict a initial velocity through the neural network and compute the most probable transition pathway $U$. If the end position of the path $U$ located in the $\epsilon$-neighborhood of the given point $(x,y)$, the transition pathway is said to have successfully transit to the limit cycle. Here, $\epsilon$ is a relatively small positive real number, representing the threshold of the limit cycle. 
 \begin{figure}[htbp]
    \centering
    \includegraphics[width=6.5cm]{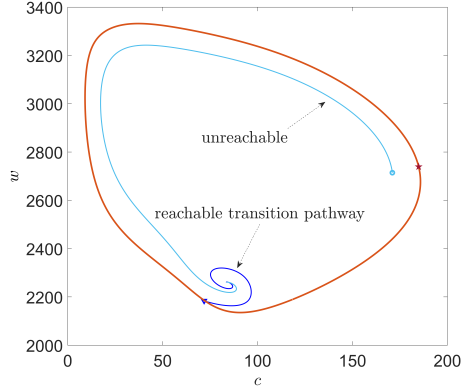}
    \caption{Display of transition pathway with parameter $\nu$=0.1. The dark blue one  successfully transits, and the light blue one fails to have a transition from stable point to limit cycle, since it does not satisfy the judgment condition.}\label{fig:3}
\end{figure}
 
To find the most probable transition pathway, we go through the following steps.\\

$\diamond$ Define the threshold  $\epsilon$ of the limit cycle $M$ and select the reachable transition pathways. A transition pathway $U$ is said to be reachable if and only if it satisfies the following condition. \\
\begin{equation}
(\left|U_{x}(\operatorname{end})-x\right|^2+\left|U_{y}(\operatorname{end})-y\right|^2)^\frac{1}{2}<\epsilon.
\end{equation}
Here, $U$ represents the most probable transition pathway computed by the neural network for a given point $(x,y)\in M$, $M$ is the limit cycle, and $U_x$, $U_y$ are two components of $U$. An example is shown in Fig. \ref{fig:3}.

$\diamond$ By calculating the Onsager-Machlup functional value of each path, the minimum functional value is selected, that is, the most probable transition pathway according to section 2.2.

\begin{figure}[htbp]

\begin{minipage}[]{0.5 \textwidth}
 \leftline{~~~~~\tiny\textbf{(\ref{fig:4}a)}}
\centerline{\includegraphics[width=6.5cm]{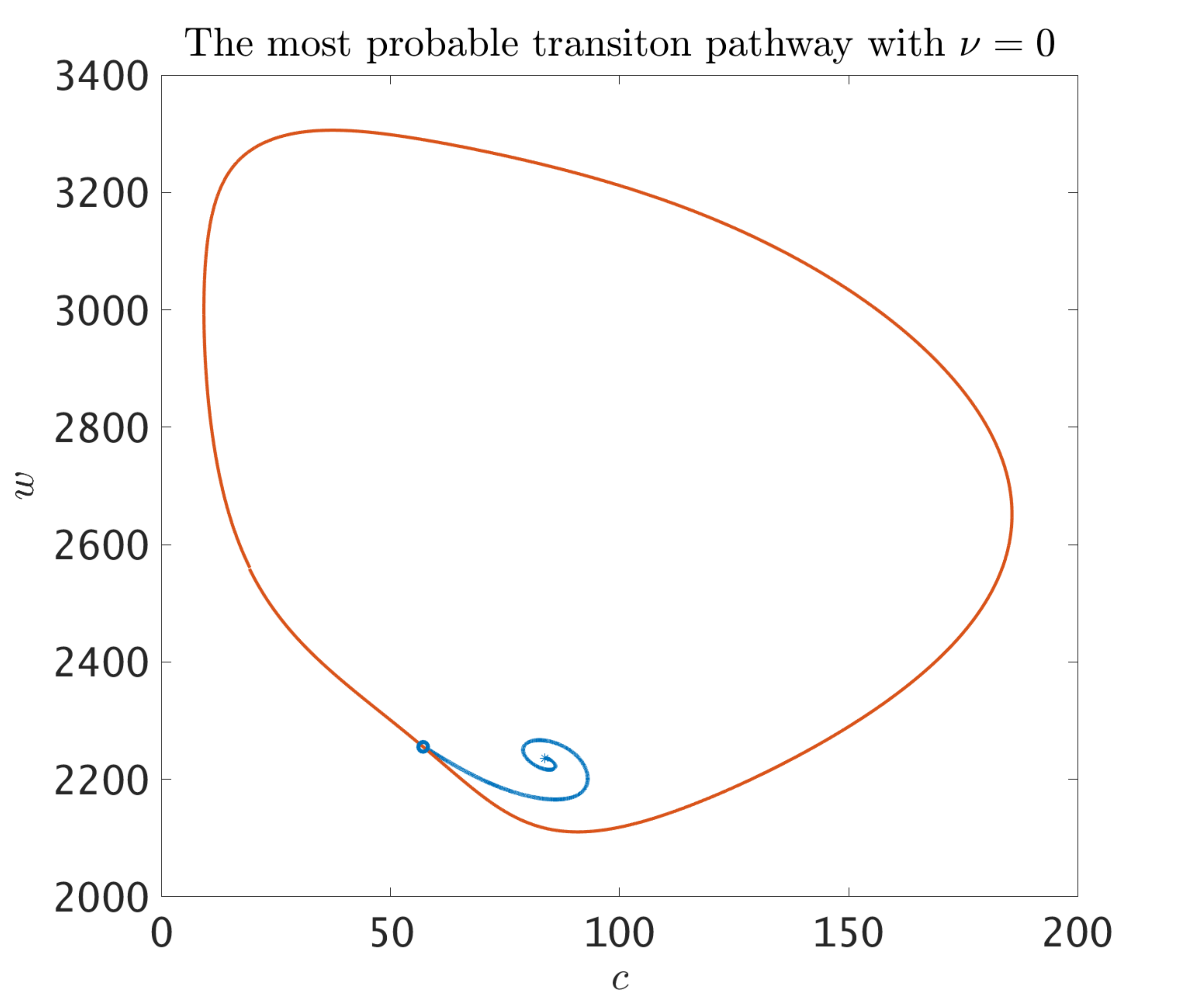}}
\end{minipage}
\hfill
\begin{minipage}[]{0.5 \textwidth}
 \leftline{~~~~~\tiny\textbf{(\ref{fig:4}b)}}
\centerline{\includegraphics[width=6.5cm]{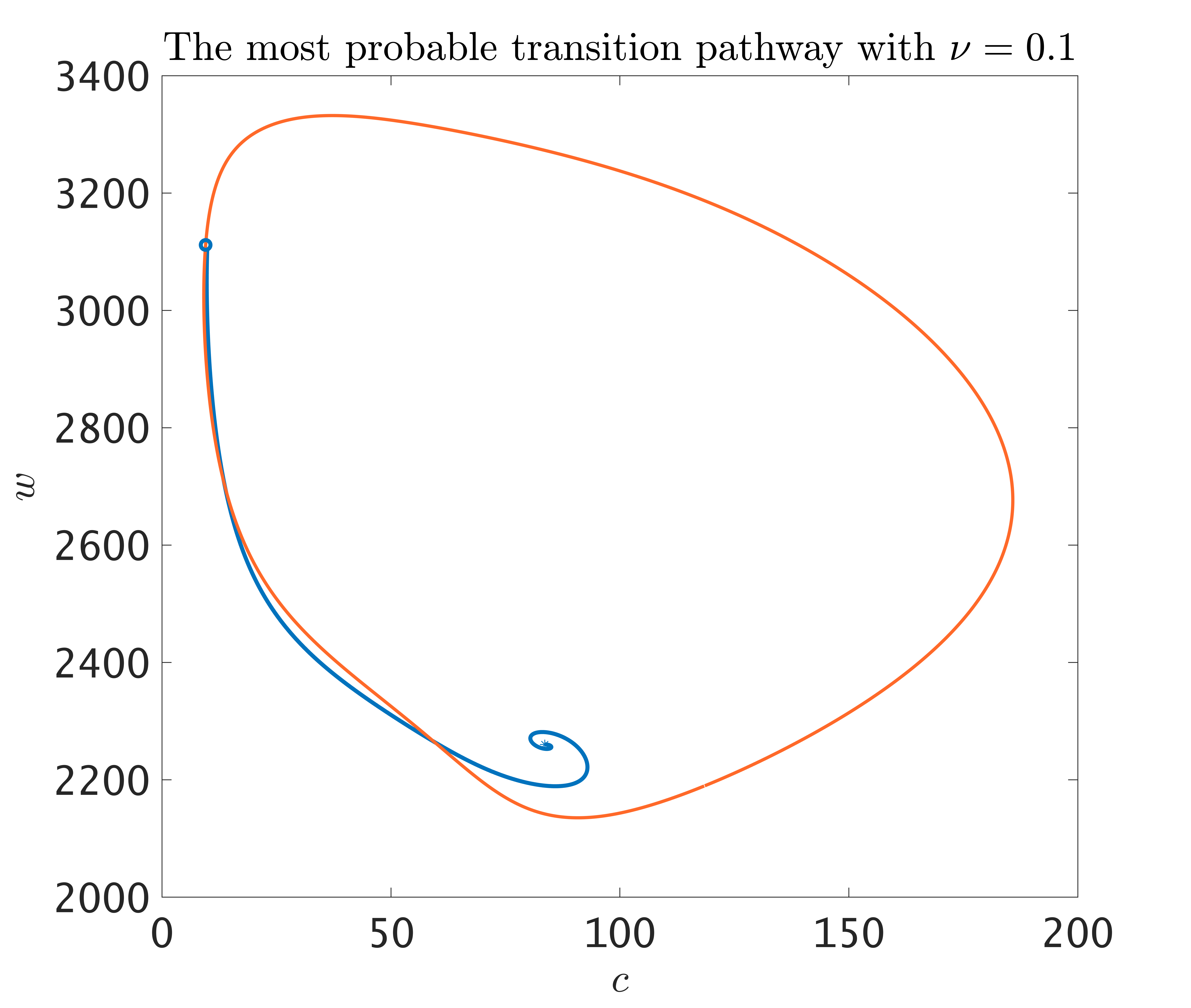}}
\end{minipage}
\begin{minipage}[]{0.5 \textwidth}
 \leftline{~~~~~\tiny\textbf{(\ref{fig:4}c)}}
\centerline{\includegraphics[width=6.5cm]{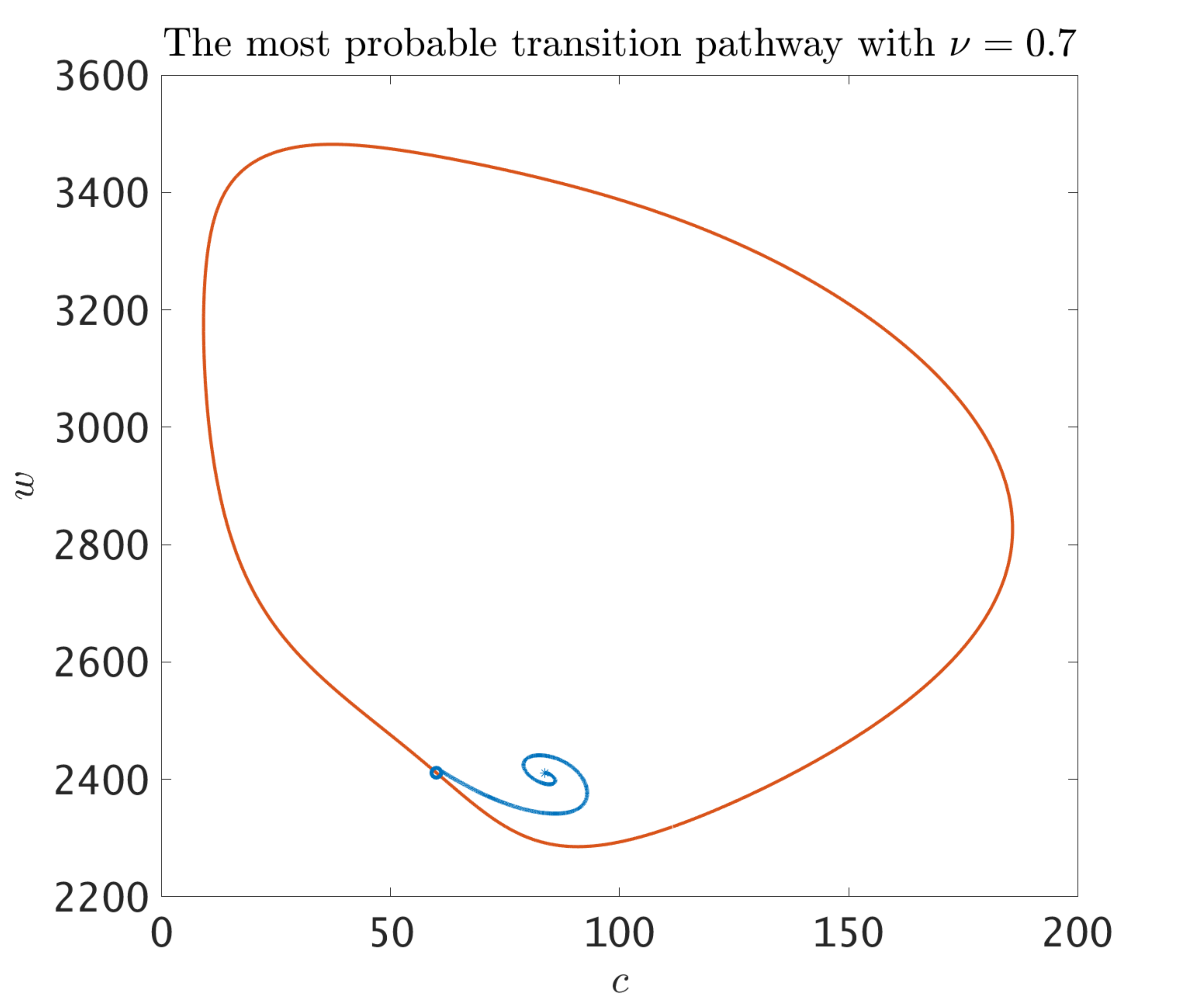}}
\end{minipage}
\hfill
\begin{minipage}[]{0.5 \textwidth}
 \leftline{~~~~~\tiny\textbf{(\ref{fig:4}d)}}
\centerline{\includegraphics[width=6.5cm]{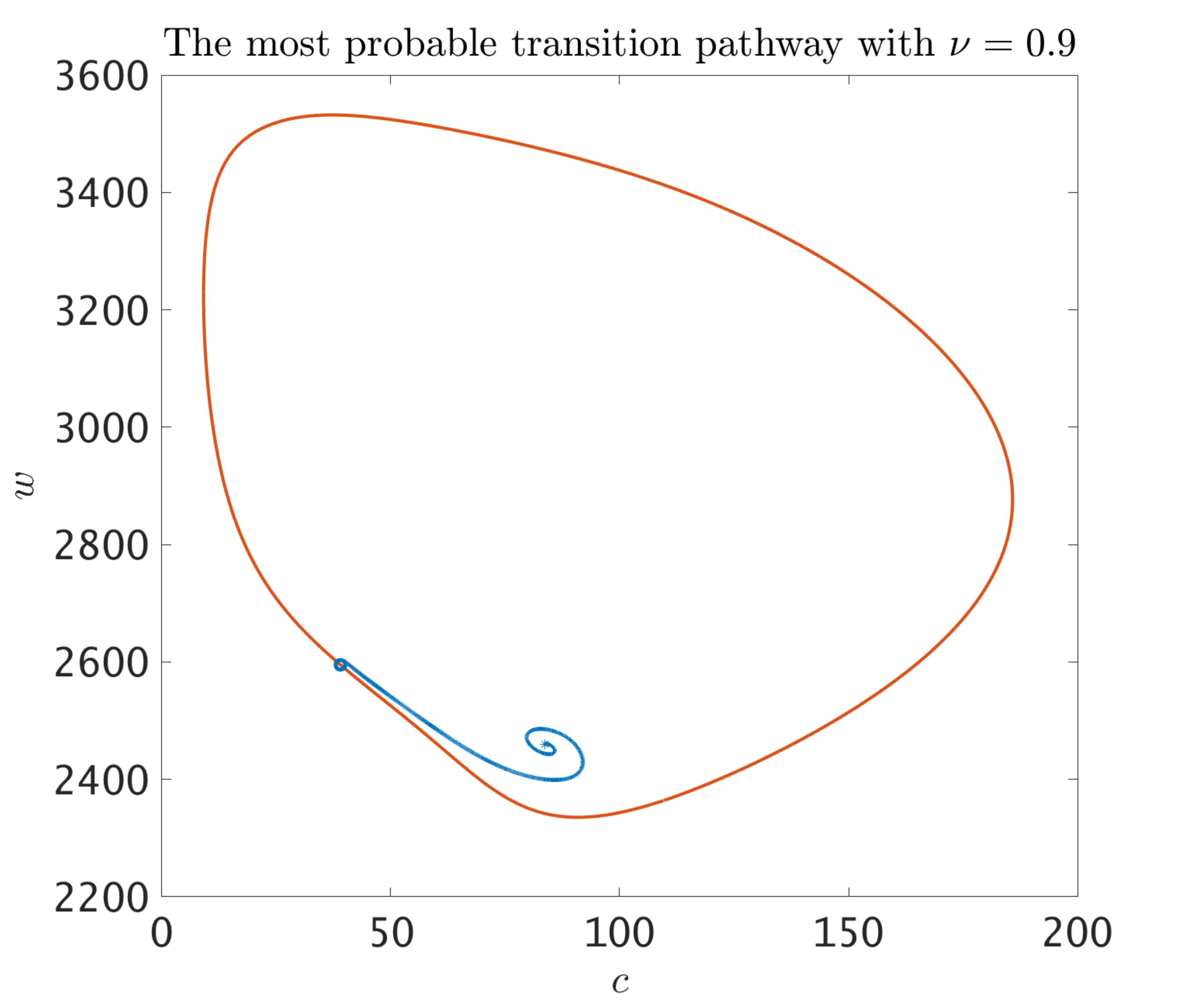}}
\end{minipage}
\caption{The most probable transition pathway with respect to different input rates $\nu$. (4a) $\nu$=0; (4b) $\nu$=0.1; (4c) $\nu$=0.7; (4d) $\nu$=0.9.}\label{fig:4}
\end{figure}

We display the most probable transition pathway under four different input rates $\nu$ as shown in Fig. \ref{fig:4}.

It shows the most probable transition pathway (solid blue line) for the carbon cycle system from a stable point (blue asterisk) to an oscillatory state (red circle) at different carbon input rates $\nu$. According to section $3$, it is the most probable transition pathway from the view of Onsager-Machlup functional, that is, the blue pathway is the one with the smallest Onsager-Machlup functional value and owns the highest probability of occurrence among all transition pathways. Fig. \ref{fig:4} is showing us a unified trend that as $\nu$ changes, the carbon cycle system tends to transit to a state which has a smaller $c$ (except for $\nu$=0.1).

\begin{figure}[htbp]
\begin{minipage}[]{0.5 \textwidth}
 \leftline{~~~~~\tiny\textbf{(\ref{fig:5}a)}}
\centerline{\includegraphics[width=6.5cm]{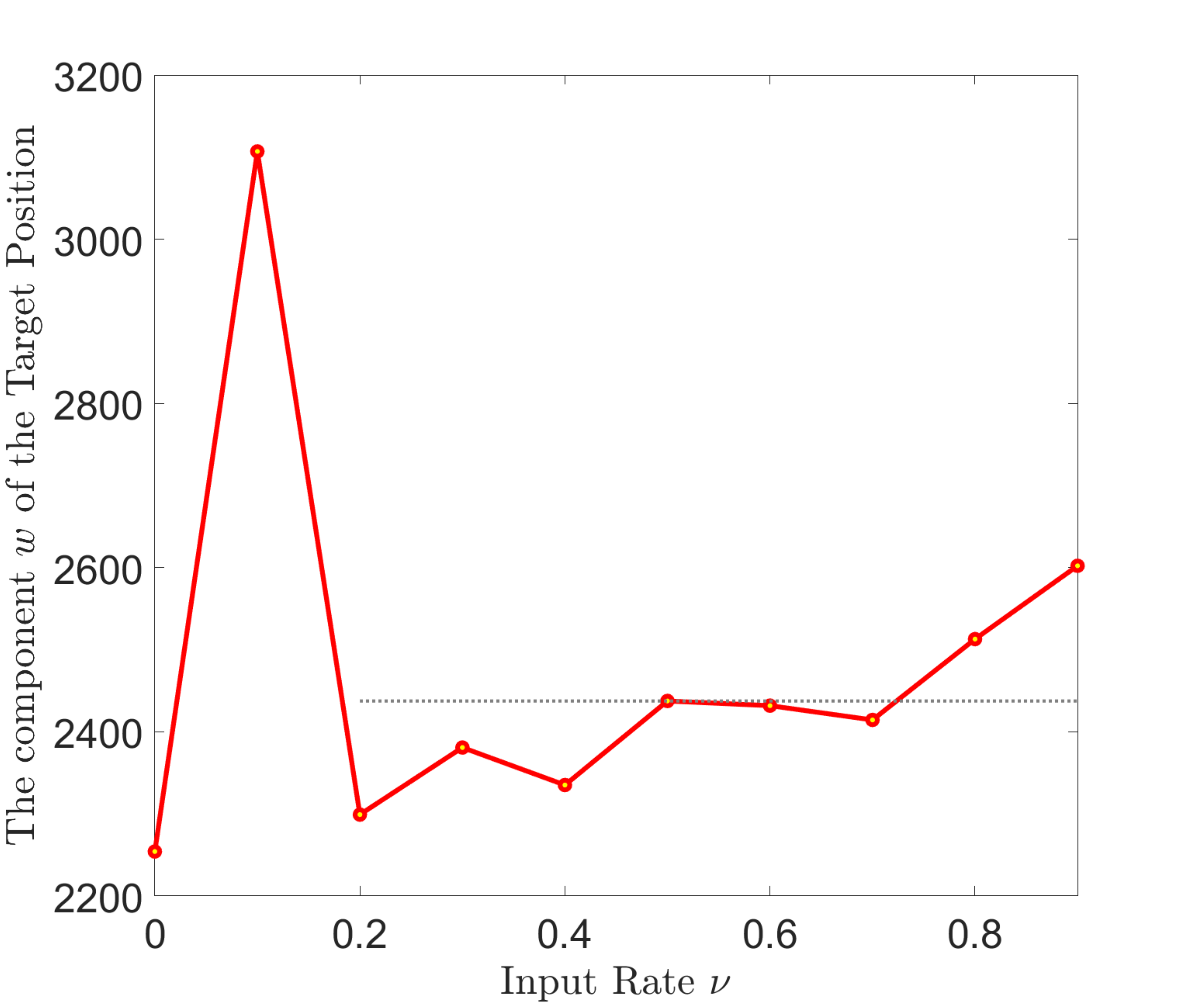}}
\end{minipage}
\hfill
\begin{minipage}[]{0.5 \textwidth}
 \leftline{~~~~~\tiny\textbf{(\ref{fig:5}b)}}
\centerline{\includegraphics[width=6.5cm]{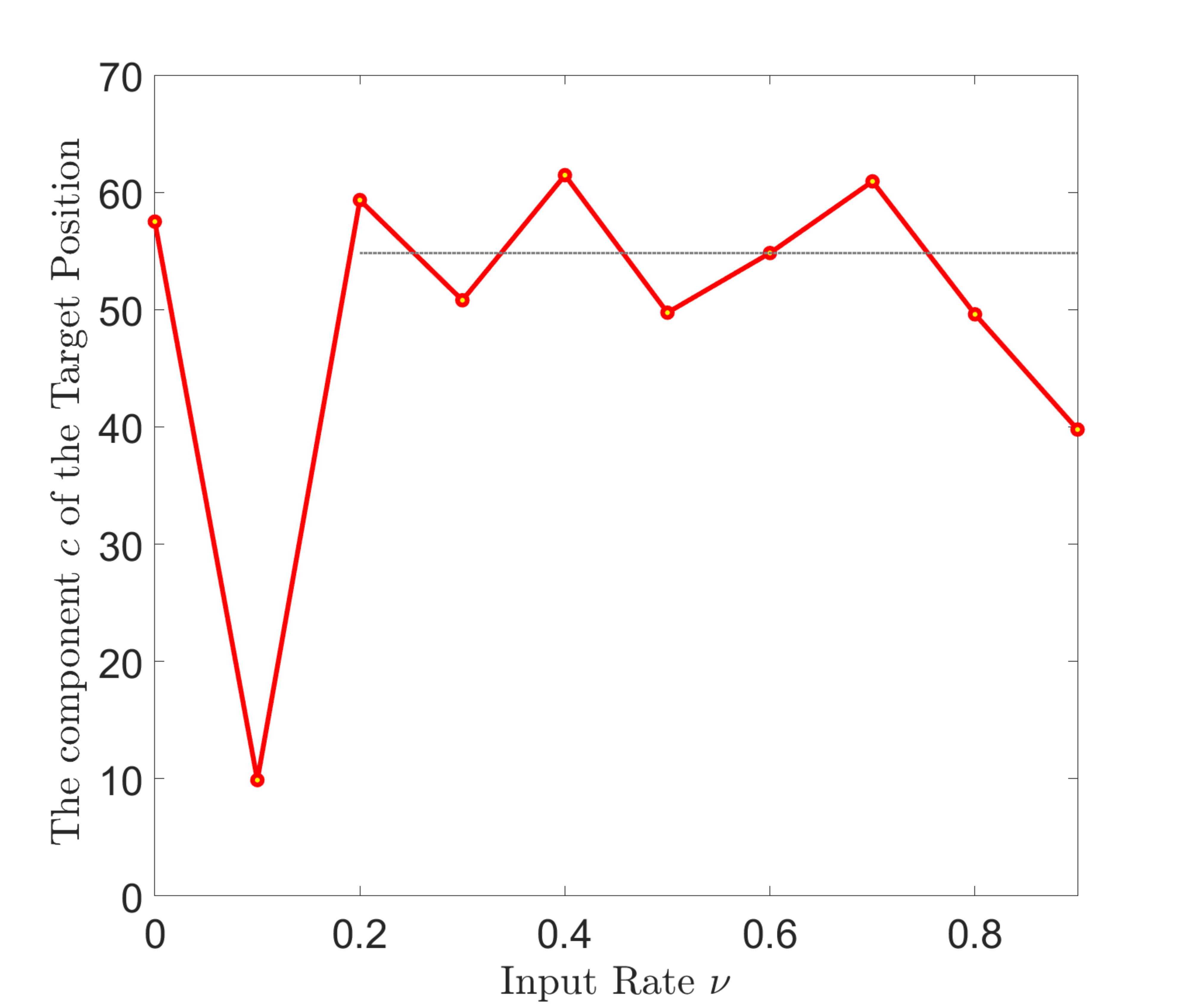}}
\end{minipage}
\caption{The components $CO_3^{2-}$ concentration $c$ and $DIC$ concentration $w$ of the target position with respect to different external input rates $\nu$. (5a): The concentration $w$ of the target position with respect to $\nu$. (5b): The $CO_3^{2-}$ concentration $c$ of the target position with respect to $\nu$.}\label{fig:5}
\end{figure}

Fig. \ref{fig:5} shows the relationship between different parameters $\nu$ and the end point components $c$ and $w$,respectively.

As shown above, not only $w$, but also $c$ is fluctuated with $\nu$, and there has been a huge oscillation at $\nu$=0.1, which can be considered as a critical transition. At this time, the concentration of $CO_3^{2-}$ falls down to a local minimal value 9.8782 $\mu mol\cdot kg^{-1}$.  When $\nu$ is growing from 0.2, the value of $CO_3^{2-}$ concentration $c$ fluctuates around $58$ $\mu mol\cdot kg^{-1}$ without a mutation.

According to section 2, the most probable transition pathway is characterized by Onsager-Machlup action functional, since the parameter $\nu$ shows such obvious effect on the transition pathway, it should also affect Onsager-Machlup action functional. We make this reasonable assumption, since there is an oscillation at $\nu$=0.1, the value of the Onsager-Machlup action functional will behave differently at $\nu$=0.1.

\begin{figure}[htbp]
    \centering
    \includegraphics[width=6.5cm]{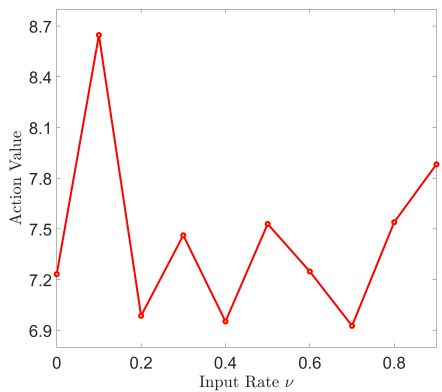}
    \caption{The action value of the most probable transition pathway with respect to different external input rates $\nu$.}\label{fig:6}
\end{figure}

Fig. \ref{fig:6} shows the relationship between the action functional value and  the most probable transition pathway at different input rate $\nu$. Because of the limit of the space, we only show four results, more results can be referred to \url{https://github.com/Cecilia-ChenJY/carbon-cycle-system}. There is no accident that at $\nu$=0.1, it experienced a big oscillation, with the value rising first, then falling, and finally fluctuating around the value 7.2 $\mu mol\cdot kg^{-1}$. And this can be considered as a critical tipping point. And as our assumption goes, the Onsager-Machlup action functional reaches the top value at $\nu$=0.1. This phenomenon indicates that when the external input rate of carbon dioxide is $\nu$=0.1, a large amount of energy is needed to transit, and the probability of state transition is lower at this time.

\subsection{Analysis of effects of transition time on the most probable transition pathway}
The reader may wonder why the neural network shows that the most probable transition pathway end up on the left side of the limit cycle, but what does the path to the right look like?  Is it possible that the path could transit to the right side of the limit cycle? We may think it possible. The transition path is not only governed by velocity, but also by the transition time $T$.

As for the neural network, the error of the model is not a good interpretation of the velocity error of a single point, because it is the sum of the errors of all points. We reasonably guess that the transition pathway whose end point falls on the right cannot reach the limit cycle due to its high dependence and sensitivity to velocity.
There are two ways to solve this problem: (i) increase the threshold of the limit cycle, that is, increase $\epsilon$, (ii) extend $T$ to see if it can be reached within large enough time. 

In this section, we detect the impacts of transition time $T$ on transition pathway. Actually, we fix the external carbon input rate $\nu$.

\begin{figure}[htbp]

\begin{minipage}[]{0.5 \textwidth}
 \leftline{~~~~~\tiny\textbf{(\ref{fig:7}a)}}
\centerline{\includegraphics[width=6.5cm]{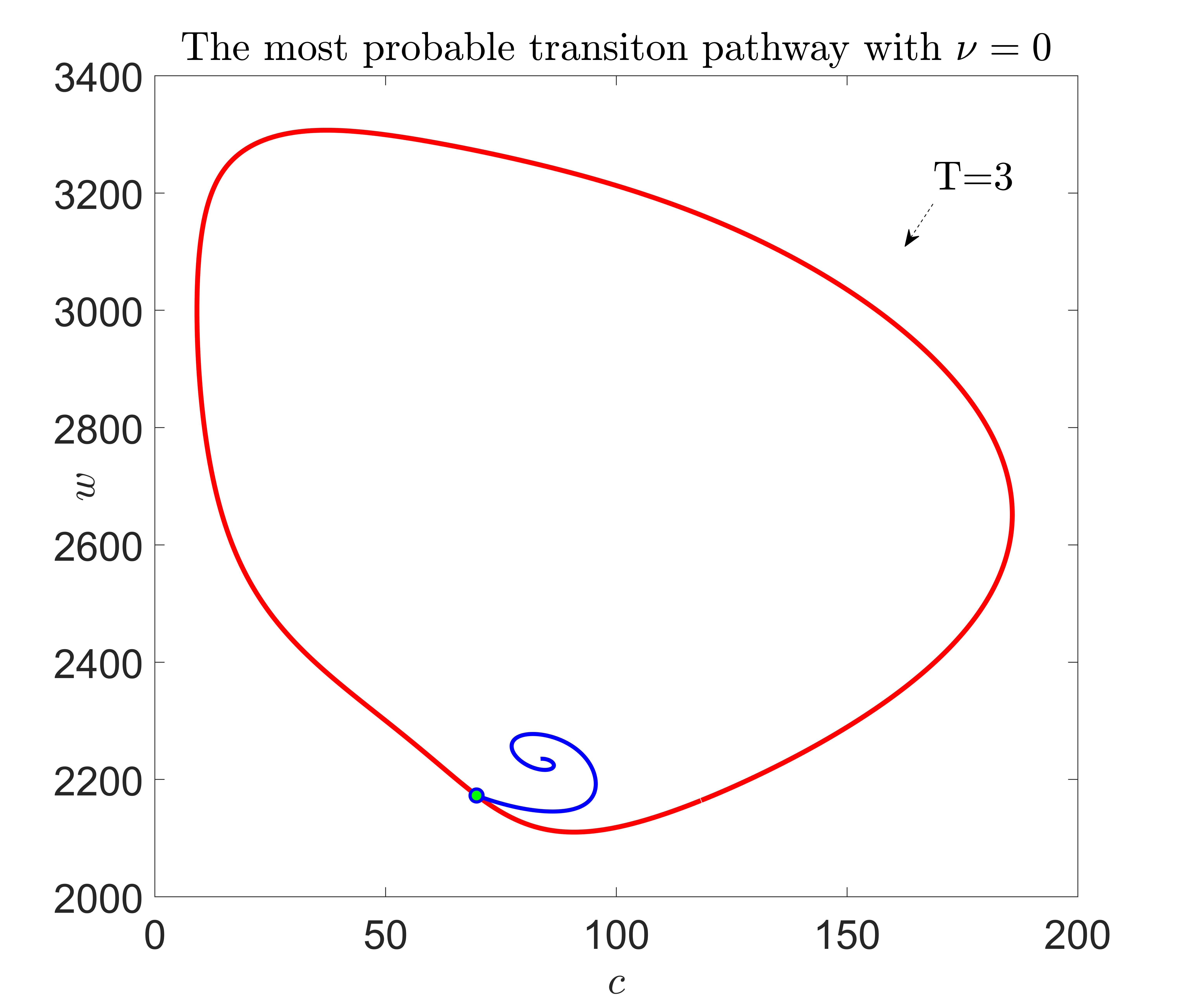}}
\end{minipage}
\hfill
\begin{minipage}[]{0.5 \textwidth}
 \leftline{~~~~~\tiny\textbf{(\ref{fig:7}b)}}
\centerline{\includegraphics[width=6.5cm]{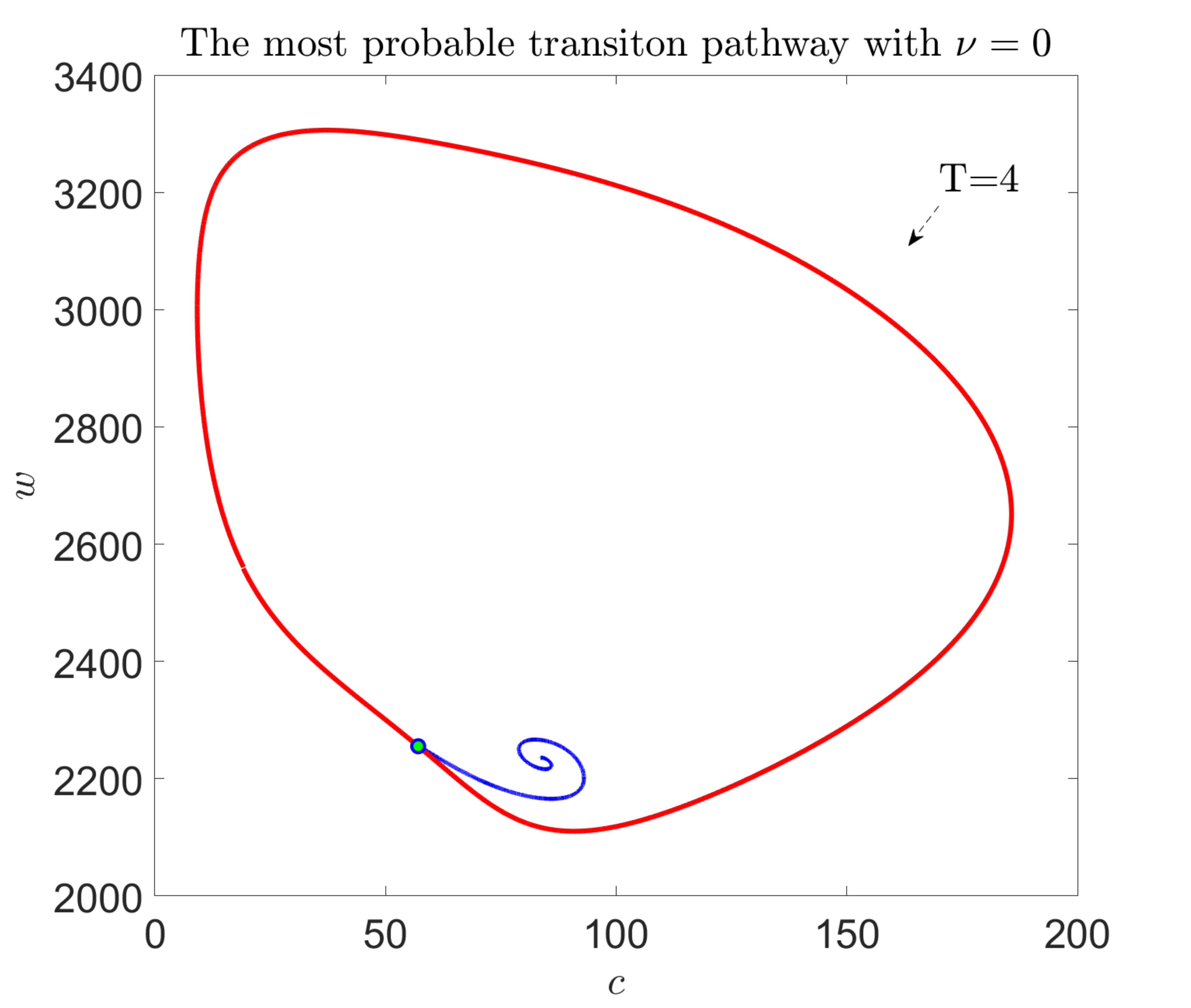}}
\end{minipage}
\begin{minipage}[]{0.5 \textwidth}
 \leftline{~~~~~\tiny\textbf{(\ref{fig:7}c)}}
\centerline{\includegraphics[width=6.5cm]{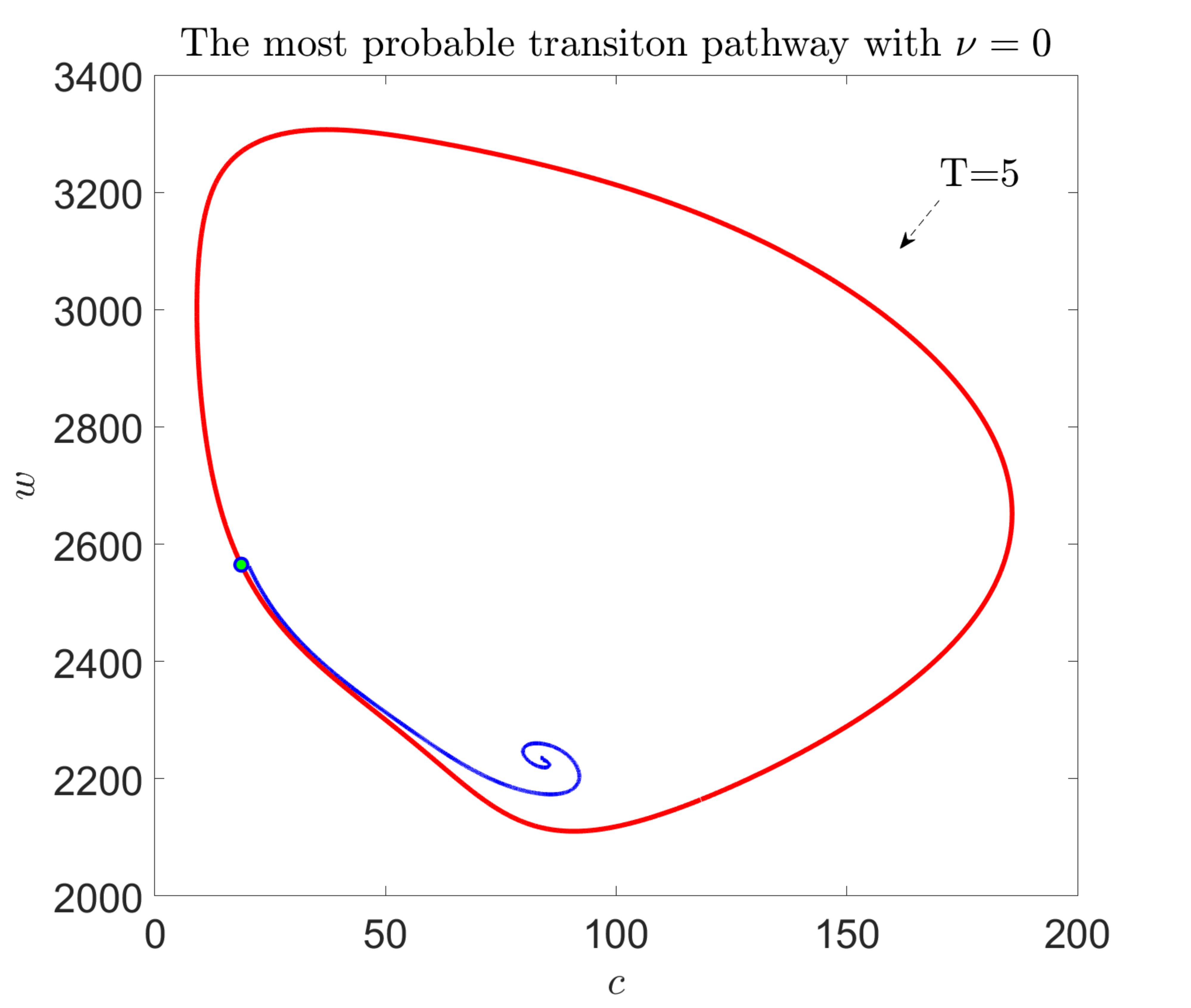}}
\end{minipage}
\hfill
\begin{minipage}[]{0.5 \textwidth}
 \leftline{~~~~~\tiny\textbf{(\ref{fig:7}d)}}
\centerline{\includegraphics[width=6.5cm]{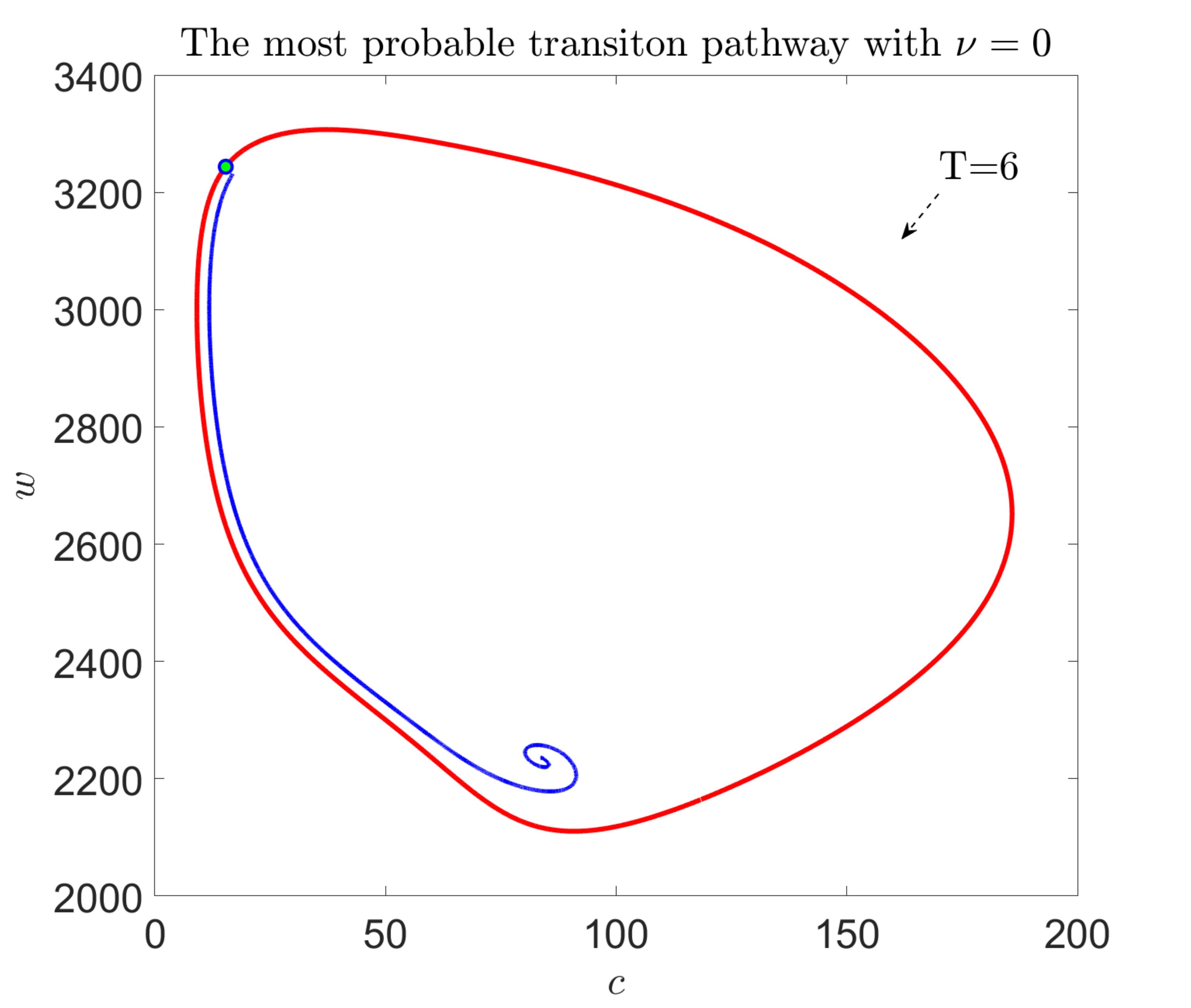}}
\end{minipage}
\caption{The impacts of time $T$ on the most probable transition pathway with $\nu$=0. (7a) $T$=3; (7b) $T$=4; (7c) $T$=5; (7d) $T$=6.}\label{fig:7}
\end{figure}

We take the simplest case where we have no external carbon dioxide input, we fix the parameter $\nu$=0. The same method is taken as in subsection $3.1$, the most probable transition pathways at different time are computed respectively shown in Fig. \ref{fig:7}. We analyze the results in Fig. \ref{fig:7}, an interesting phenomenon indicates that as time increases, the endpoint of the most probable transition pathway moves toward the right half of the limit cycle. That is, for the stochastic carbon cycle system, with the increase of time, the system is more and more inclined to transition from a stable state to an oscillatory state with increasing $CO_3^{2-}$ concentration $c$. Concentration of carbon dioxide $c$ is going to go down and then up, and we can make an assumption, as long as $T$ is big enough, the most probable transition pathway is going to go to the right half of the limit cycle. And a recently result from \cite{wei2021most} confirms that.

\subsection{Forecast of the most probable transition time}
Since time does affect the transition of the carbon cycle system, people also care about the time it takes for a transition event. It is also important to forecast the transition time when a disruption comes with the largest probability. We shall focus on the most violent change of $CO_3^{2-}$ concentration in the limit cycle, which locates at `1722' of the limit cycle. To forecast the optimal transition time, we divide up the transition time $T$ in [1,11] into 200 pieces and use the PINNs method to compute the most probable transition pathway. Then, the optimal transition time is the one with the least amount of time occupied by the minimal value of the action functional.

The PINNs method is successfully applied to compute the most transition pathway, which is an effective method, even for high dimensional systems. The idea is to consider the governing equation as a part of the loss function. 

We approximate the solution to the Euler-Lagrange equation \eqref{EL_eqn} from a set of training data, which consist of residual data, initial and final data. 

\begin{figure}[htbp]
\begin{minipage}[]{0.5 \textwidth}
  \leftline{~~~~~\tiny\textbf{(\ref{fig:8}a)}}
\centerline{\includegraphics[width=6.5cm]{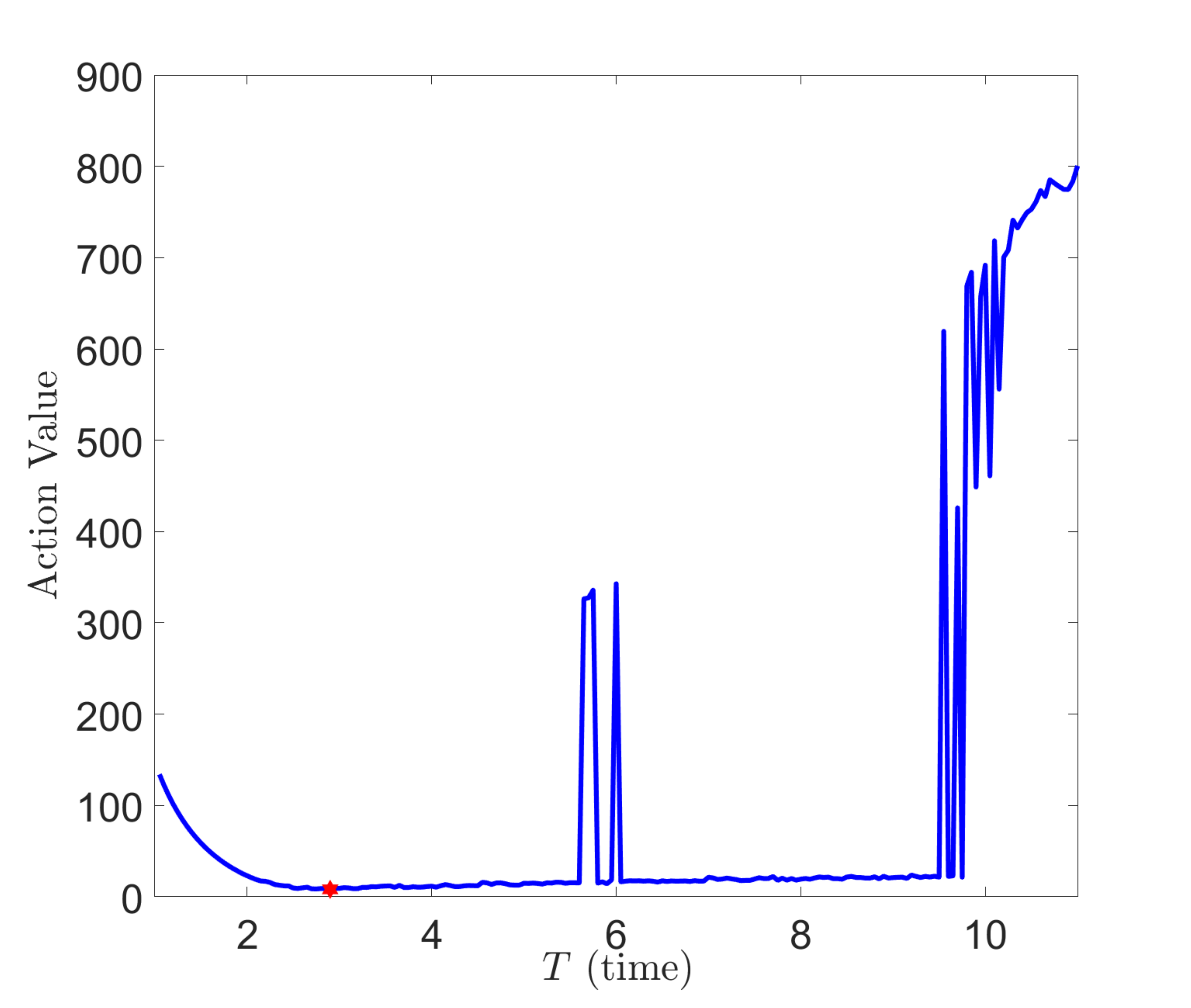}}
\end{minipage}
\hfill
\begin{minipage}[]{0.5 \textwidth}
  \leftline{~~~~~\tiny\textbf{(\ref{fig:8}b)}}
\centerline{\includegraphics[width=6.5cm]{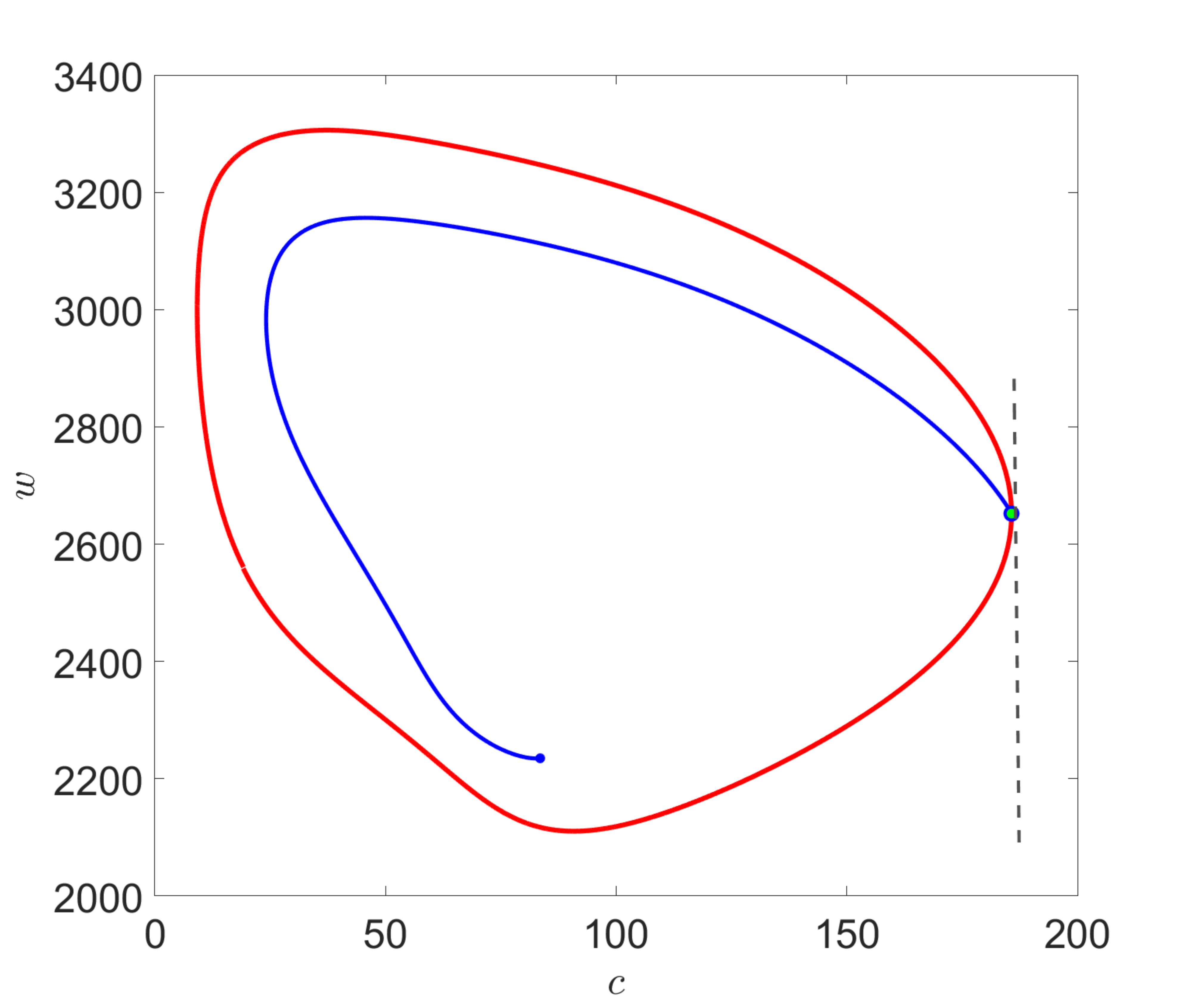}}
\end{minipage}
\caption{The action functional value and the most probable transition pathway with $\nu$=0. (8a) The variation of action functional value with time, and red point $*$ corresponds to the optimal transition time 2.9 (29 thousand years) corresponding to the minimum value. (8b) The most probable transition pathway for the optimal transition time 2.9, corresponding to the target location `1722'(the vertical dashed line) where the violent $CO_3^{2-}$ concentration change takes.}
\label{fig:8}
\end{figure}

Denote the residual data as $\{(t_j,\frac{d}{dt}\partialx{\dot{z}}L(z(t_j),\dot{z}(t_j))-\partialx{z}L(z(t_j),\dot{z}(t_j))\}_{j=1}^m$, where $t_j\in U=(0,T)$, and the initial and final data $\{(0,x_0),(T,x_T)\}$, where $m$ is the number of the residual data. Therefore the the empirical PINN loss is
\begin{equation}
\begin{split}\label{emp_loss}
\operatorname{Loss}_{m}^{\operatorname{PINN}}\left(z,\lambda\right)&=\frac{1}{m}\sum\limits_{j=1}^m\|\frac{d}{dt}\partialx{\dot{z}}L(z(t_j),\dot{z}(t_j))-\partialx{z}L(z(t_j),\dot{z}(t_j))\|^2\\
&\ \ + \frac{\lambda}{2}(\|h(0)-x_0\|^2+\|h(T)-x_T\|^2),
\end{split}
\end{equation}
where $\lambda$ is a positive constant, which is used to balance the residual loss and boundary loss. 

We shall use a full connected neural network to compute the most transition pathway by optimizing the empirical PINN loss \eqref{emp_loss} with $\lambda = 60$. All the neural networks have 4 hidden layers and 20 neurons per layer, with $\operatorname{tanh}$ activation function. The weights are initialized with truncated normal distributions and the biases are initialized as zero.
The Adam optimizer with a learning rate of $10^{-3}$ is used to train the loss function. In one dimension, the number of residual points for evaluating the Euler-Lagrange equation is $m=501$. 

The results are shown in Fig. \ref{fig:8}. It presents the action functional value respect to the transition time. When the transition time $T=2.9$, it reaches the minimum value, which corresponds to the most probable occurrence of the transition. We obtain the optimal transition time $T=2.9$ (29 thousand years). It provides us an early warning to forecast the transition. The most transition pathway is shown in Fig. (\ref{fig:8}b).

\section{Conclusion}
Carbon cycle system is an important indicator for revealing climate changes. It reacts to the occurrence of emergencies, such as volcanic eruption and human activities. In order to illustrate how this works in our climate system, we examine a stochastic dynamical system model of the marine carbon cycle. We devote to the study of the most probable transition pathway between a metastable state and the oscillatory regimes in this stochastic carbon cycle system, and we carry out the numerical calculation by a neural shooting and a physics informed neural network methods. The results show that when the external carbon input rate $\nu$ of the carbon cycle system is disturbed by a certain intense of noise, different $\nu$ corresponds to the different transition pathway. However, things are different at $\nu$=0.1.

Further, we discuss the effects of transition time on the most probable transition pathway, and the results show that when the external carbon input rate $\nu$ is fixed, the most probable transition pathway of the stochastic carbon cycle system tends to shift towards the state of where $CO_3^{2-}$ concentration $c$ becoming larger over time. Which may be related to the natural environment, and the evolution of time brings a higher concentration of $CO_3^{2-}$ $c$. In order to forecast the occurrence of some events (usually catastrophic emergencies) and provide the early warning, we use a physics informed neural network to solve the optimal transition time at a fixed point in the state of oscillatory regimes, which is of great significance to the prediction of catastrophic events.

These results imply that this stochastic carbon cycle system can detect the phenomenon of transition between a metastable state and oscillatory regimes. For random perturbations are ubiquitous in nature, our study has strong practical consequences. In addition, the transition time has an effect on the variability of $CO_3^{2-}$ concentration $c$, which motivates us to control carbon emissions in the environment. Even more useful is our contribution to the prediction of unexpected events, as this paper calculates the optimal time at which state with the most violent change in  $CO_3^{2-}$ concentration $c$, providing early warning to prevent catastrophic events from sudden changes in the environment.

\section*{Acknowledgements}
We would like to thank Haokun Li, Ting Gao for helpful discussions. This work was partly supported by NSFC grants 11771449.

\section*{Data Availability}
The data that support the findings of this study are openly available in GitHub.\\
\url{https://github.com/Cecilia-ChenJY/carbon-cycle-system}

\bibliographystyle{abbrv}
\bibliography{References}
\begin{appendices}
\section{Derivation of Equation (\ref{EL_eqn})}\label{AppA}
For the stochastic carbon cycle system \eqref{carbon}, we shall derive the Euler-Lagrangian equation \eqref{EL_eqn} in detail. Denote the diffusion matrix $\sigma(x,y)=
\left[
\begin{matrix}
   -\mu f(x)  &  0\\
   0  & \mu
\end{matrix}
\right]
$ and the vector field $\tilde{b}=(\tilde{b}^1,\tilde{b}^2)$, where 
\begin{align*}
\tilde{b}^1(x,y)&=f(x)\left\{\mu\left[1-b s\left(x, c_{p}\right)-\theta \bar{s}\left(x, c_{x}\right)-\nu \right]+y-w_{0}\right\},\\
\tilde{b}^2(x,y)&=\mu\left[1-b s\left(x, c_{p}\right)+\theta \bar{s}\left(x, c_{x}\right)+\nu\right]-y+w_{0}.
\end{align*}
Thus, the Riemanian metric is $V(x,y)=(g_{ij})_{2\times2}(x,y)=
\left[
\begin{matrix}
   \frac{1}{\mu^2f^2(x)}  &  0\\
   0  & \frac{1}{\mu^2}
\end{matrix}
\right]
$ and its inverse matrix is $(g^{ij})_{2\times2}(z)=
\left[
\begin{matrix}
   {\mu^2f^2(x)}  &  0\\
   0  & {\mu^2}
\end{matrix}
\right]
$ . The determinate of $V$ is $|V(x,y)|=\frac{1}{\mu^4f^2(x)}$.

\noindent\textbf{Compute the Christoffel symbols}
\begin{align*}
 \Gamma_{11}^1&=\frac{1}{2}g^{11}\left(\partialx{x}g_{11}+\partialx{x}g_{11}-\partialx{x}g_{11}\right)+\frac{1}{2}g^{12}\left(\partialx{x}g_{12}+\partialx{x}g_{21}-\partialx{x}g_{11}\right)\\
  &=\frac{1}{2}\mu^2 f^2(x)\cdot\partialx{x}\left(\frac{1}{\mu^2 f^2(x)}\right)=\frac{1}{2} f^2(x)\cdot\partialx{x}\left(\frac{1}{ f^2(x)}\right)\\
  &=-\frac{f'(x)}{f(x)},\\
 \Gamma_{12}^1&=\frac{1}{2}g^{11}\left(\partialx{x}g_{21}+\partialx{y}g_{11}-\partialx{x}g_{12}\right)+\frac{1}{2}g^{12}\left(\partialx{x}g_{22}+\partialx{y}g_{21}-\partialx{y}g_{12}\right)\\
  &=0=\Gamma_{21}^1,\\
  \Gamma_{22}^1&=\frac{1}{2}g^{11}\left(\partialx{y}g_{21}+\partialx{y}g_{12}-\partialx{x}g_{22}\right)+\frac{1}{2}g^{12}\left(\partialx{y}g_{22}+\partialx{y}g_{22}-\partialx{y}g_{22}\right)\\
  &=0,\\
  \Gamma_{11}^2&=\frac{1}{2}g^{21}\left(\partialx{x}g_{11}+\partialx{x}g_{11}-\partialx{x}g_{11}\right)+\frac{1}{2}g^{22}\left(\partialx{x}g_{12}+\partialx{x}g_{21}-\partialx{y}g_{11}\right)\\
  &=-\frac{1}{2} \mu^2\cdot\partialx{y}\left(\frac{1}{\mu^2\cdot f^2(x)}\right)\\
  &=0,\\
  \Gamma_{12}^2&=\frac{1}{2}g^{21}\left(\partialx{x}g_{21}+\partialx{y}g_{11}-\partialx{x}g_{12}\right)+\frac{1}{2}g^{22}\left(\partialx{x}g_{22}+\partialx{y}g_{21}-\partialx{y}g_{12}\right)\\
  &=0=\Gamma_{21}^2,\\
  \Gamma_{22}^2&=\frac{1}{2}g^{21}\left(\partialx{x}g_{21}+\partialx{y}g_{11}-\partialx{x}g_{12}\right)+\frac{1}{2}g^{22}\left(\partialx{y}g_{22}+\partialx{y}g_{22}-\partialx{y}g_{22}\right)\\
  &=0.
\end{align*}
\noindent\textbf{Compute the divergence}

The modified vector field is given by 
\begin{align*}
    b^1&=\tilde{b}^1-\frac{1}{2}\sum_{lj}(V^{-1}(z))^{lj} \Gamma^1_{lj}=\tilde{b}^1-\frac{1}{2}\mu^2f^2(x)\cdot\left(-\frac{f'(x)}{f(x)}\right)=\tilde{b}^1+\frac{1}{2}\mu^2f'(x)f(x),\\
    b^2&=\tilde{b}^2-\frac{1}{2}\sum_{lj}(V^{-1}(z))^{lj} \Gamma^2_{lj}=\tilde{b}^2.
\end{align*}
Therefore, the divergence of the modified vector field is
\begin{align*}
   \mathtt{div}b(z)&=\frac{1}{\sqrt{|V(z)|}}\sum_{i}\partialx{z_i}(b^i(z)\sqrt{|V(z)|})\\
   &=\mu^2|f(x)|\cdot \left(\partialx{x}\left(b^1(z)\cdot\frac{1}{\mu^2|f(x)|}\right)+
   \partialx{y}\left(b^2(z)\cdot\frac{1}{\mu^2|f(x)|}\right)\right)\\
   &=\partialx{x}b^1(z)+\partialx{y}b^2(z)+b^1(z)|f(x)|\cdot\partialx{x}\frac{1}{|f(x)|}\\
   &=\partialx{x}b^1(z)+\partialx{y}b^2(z)-b^1(z)\cdot\left(\frac{f'(x)}{f(x)}\right).
\end{align*}
\noindent\textbf{Compute the scalar curvature}

We use the notation of Einstein sum in the following computing. The scalar curvature is defined as the trace of the Ricci curvature tensor with respect to the Riemanian metric
\begin{align*}
R=g^{ij}R_{ij},
\end{align*}
where $R_{ij}$ is the Ricci curvature tensor. By \cite[Eqn 4.1.32, Eqn 4.3.6, Eqn 4.3.18]{jost2008riemannian}, we have
\begin{align*}
R&=g^{ik} R_{ik}=g^{ik}(g^{jl} R_{ijkl})=g^{ik}(g^{jl} (g_{im}R^{m}_{jkl}))\\
&=g^{ik}(g^{jl}(g_{im}(\partialx{x^k}\Gamma_{lj}^m-\partialx{x^l}\Gamma_{kj}^m+\Gamma_{ka}^m\Gamma_{lj}^a-\Gamma_{la}^m\Gamma_{kj}^a)).
\end{align*}
Since $g^{jl}g_{im}=\delta_{i=j,l=m}$, we have
\begin{align*}
R&=g^{ik}R^{m}_{ikm}=g^{ik}(\partialx{x^k}\Gamma_{im}^m-\partialx{x^m}\Gamma_{ki}^m+\Gamma_{ka}^m\Gamma_{mi}^a-\Gamma_{ma}^m\Gamma_{ki}^a).
\end{align*}
Since $\Gamma_{11}^1=-\frac{f'(x)}{f(x)}$, $\Gamma_{jk}^i=0 (\text{others})$, $g^{11}=\mu^{2}f^{2}(x)$, $g^{22}=\mu^{2}$, and $g^{12}=g^{21}=0$, we have
\begin{align*}
R=&g^{11}(\partialx{x}\Gamma_{1m}^m-\partialx{x^m}\Gamma_{11}^m+\Gamma_{1a}^m\Gamma_{m1}^a-\Gamma_{ma}^m\Gamma_{11}^a)\\
&+g^{22}(\partialx{y}\Gamma_{2m}^m-\partialx{x^m}\Gamma_{22}^m+\Gamma_{2a}^m\Gamma_{m2}^a-\Gamma_{ma}^m\Gamma_{22}^a)\\
=&g^{11}(\partialx{x}\Gamma_{11}^1-\partialx{x}\Gamma_{11}^1+\Gamma_{11}^1\Gamma_{11}^1-\Gamma_{11}^1\Gamma_{11}^1)\\
=&0.
\end{align*}
\noindent\textbf{Compute the Euler-Lagrange equation}

The Euler-Lagrange equation in two dimension is 
\[
\frac{d}{dt}\partialx{\dot{z}}L(z,\dot{z})=\partialx{z}L(z,\dot{z})
\Rightarrow
\left\{
\begin{array}{l}
   \frac{d}{dt}\partialx{\dot{x}}L(z,\dot{z})=\partialx{x}L(z,\dot{z}),\\
   \\
    \frac{d}{dt}\partialx{\dot{z}}L(z,\dot{y})=\partialx{y}L(z,\dot{z}).
\end{array}
\right.
\]

By previous computing, the Lagrangian \eqref{EL_eqn} for the stochastic carbon cycle system \eqref{carbon} reduces to 
\begin{align*}
    L(z, \dot{z})=&(\dot{z}-b(z))V(z)(\dot{z}-b(z))+\mathtt{div}b(z)-\frac{1}{6}R(z)\\
                =&(\dot{x}-\tilde{b}^1-\frac{1}{2}\mu^2f'(x)f(x))^2\cdot \frac{1}{\mu^2f^2(x)}+\left(\dot{y}-\tilde{b}^2\right)^2\cdot\frac{1}{\mu^2} \\
                &+\partialx{x}b^1(z)+\partialx{y}b^2(z)-b^1(z)\cdot\left(\frac{f'(x)}{f(x)}\right).
\end{align*}
For simplicity, denote the  Lagrangian as 
$$
L(z,\dot{z})=a(x)(\dot{x}-b^1(x,y))^2+\frac{1}{\mu^2}(\dot{y}-b^2(x,y))^2+e(x,y),
$$
where
\[
\left\{
\begin{array}{ll}
    &a(x)=\frac{1}{\mu^2}\cdot\frac{1}{f^2(x)}, \\
    & e(x,y)=\partialx{x}b^1(z)+\partialx{y}b^2(z)-b^1(z)\cdot\left(\frac{f'(x)}{f(x)}\right).
\end{array}
\right.
\]
By calculating, we have
\begin{enumerate}[label=\circled{\arabic*}]
    \item  
\begin{align*}
    \partialxy{L}{\dot{x}}=&2a(x)(\dot{x}-b^1(x,y))\\
    \Rightarrow \frac{d}{dt}\partialxy{L}{\dot{x}}=&2\frac{\partial}{\partial x}a(x) \dot{x}(\dot{x}-b^1(x,y))+2a(x)\left(\ddot{x}-\partialxy{b^1(x,y)}{\dot{x}} \dot{x}-\partialxy{b^1(x,y)}{\dot{y}} \dot{y}\right)\\
    \partialxy{L}{x}=&\frac{\partial}{\partial x}a(x)(\dot{x}-b^1(x,y))^2-2a(x)\left(\dot{x}-b^1(x,y)\right) \frac{\partial}{\partial x}b^1(x,y)\\
    &-\frac{2}{\mu^2}\left(\dot{y}-b^2(x,y)\right)\partialx{x}b^2(x,y)+\partialx{x}e(x,y)
\end{align*}
\item 
\begin{align*}
\partialxy{L}{\dot{y}}=&\frac{2}{\mu^2}(\dot{y}-b^2(x,y))\\
\Rightarrow \frac{d}{dt}\partialxy{L}{\dot{y}}=&\frac{2}{\mu^2}\left(\ddot{y}-\partialxy{b^2(x,y)}{x}\dot{x}-\partialxy{b^2(x,y)}{y}\dot{y}\right)\\
\partialxy{L}{y}=&-2a(x)\left(\dot{x}-b^1(x,y)\right)\partialx{y}b^1(x,y)-\frac{2}{\mu^2}\left(\dot{y}-b^2(x,y)\right)\partialx{y}b^2(x,y)+\partialx{y}e(x,y)
\end{align*}

\end{enumerate}
Therefore, the Euler-Lagrangian equation for the stochastic carbon cycle system \eqref{carbon} reads
\begin{align*}
\ddot{x}=&\partialxy{b^1(x,y)}{x}\dot{x}+\partialxy{b^1(x,y)}{y}\dot{y}-\frac{a'(x)}{a(x)}\dot{x}(\dot{x}-b^1(x,y))+\frac{a'(x)}{2a(x)}(\dot{x}-b^1(x,y))^2\\
&-\left(\dot{x}-b^1(x,y)\right) \frac{\partial}{\partial x}b^1(x,y)-\frac{1}{\mu^2a(x)}\left(\dot{y}-b^2(x,y)\right)\partialx{x}b^2(x,y)+\frac{1}{2a(x)}\partialx{x}e(x,y),\\
\ddot{y}=&\partialxy{b^2(x,y)}{x}\dot{x}+\partialxy{b^2(x,y)}{y}\dot{y}-\mu^2a(x)\left(\dot{x}-b^1(x,y)\right)\partialx{y}b^1(x,y)\\
&-\left(\dot{y}-b^2(x,y)\right)\partialx{y}b^2(x,y)+\frac{\mu^2}{2}\partialx{y}e(x,y).
\end{align*}


%


\end{appendices}
\end{document}